\documentclass[10pt]{elsarticle}
\usepackage[cp1251]{inputenc}
\usepackage[english]{babel}
\usepackage{amsmath}
\usepackage{amssymb}
\usepackage{amsfonts}
\usepackage{rotating}
\usepackage{graphicx}
\usepackage{floatflt,epsfig}
\usepackage{lineno,hyperref}
\usepackage{enumerate}
\usepackage{colortbl}
\usepackage{array,tabularx,tabulary,booktabs}
\usepackage{longtable}
\usepackage{multirow}
\usepackage{wrapfig}
\usepackage{subcaption}
\usepackage[table]{xcolor}
\newcolumntype{$}{>{\global\let\currentrowstyle\relax}}
\newcolumntype{^}{>{\currentrowstyle}}

\journal{arXiv}
\newtheorem{lemma}{Lemma}

\newtheorem{theorem}{Theorem}

\newcommand{\proof}{\medskip\noindent{\bf Proof.~}}
\begin{document}
\renewcommand{\abstractname}{Abstract}
\renewcommand{\refname}{References}
\renewcommand{\tablename}{Figure.}
\renewcommand{\arraystretch}{0.9}
\thispagestyle{empty}
\sloppy

\begin{frontmatter}
\title{On eigenfunctions and maximal cliques of Paley graphs of square order\tnoteref{grant}}
\tnotetext[grant]{The reported study was funded by RFBR according to the research project 17-51-560008.
The first author is partially supported by the grant NSFC 11671258. The first and the third authors
are partially supported by RFBR according to the research project 16-31-00316.
}

\author[01,03,04]{Sergey~Goryainov\corref{cor1}}
\cortext[cor1]{Corresponding author}
\ead{44g@mail.ru}

\author[03]{Vladislav~V.~Kabanov}
\ead{vvk@imm.uran.ru}

\author[03,04]{Leonid~Shalaginov}
\ead{44sh@mail.ru}

\author[02]{Alexandr~Valyuzhenich}
\ead{graphkiper@mail.ru}

\address[01]{Shanghai Jiao Tong University, 800 Dongchuan RD. Minhang District, Shanghai, China}
\address[02]{Sobolev Institute of Mathematics, Ak. Koptyug av. 4, Novosibirsk, 630090, Russia}
\address[03]{Krasovskii Institute of Mathematics and Mechanics, S. Kovalevskaja st. 16, Yekaterinburg, 620219, Russia}
\address[04]{Chelyabinsk State University, Brat'ev Kashirinyh st. 129, Chelyabinsk,  454021, Russia}

\begin{abstract}
In this paper we find new maximal cliques of size $\frac{q+1}{2}$ or $\frac{q+3}{2}$, accordingly as $q\equiv 1(4)$ or $q\equiv 3(4)$,
in Paley graphs of order $q^2$, where $q$ is an odd prime power.
After that we use new cliques to define a family of eigenfunctions corresponding to both non-principal eigenvalues and having the cardinality of support $q+1$, which is the minimum by the weight-distribution bound.

\end{abstract}

\begin{keyword}
Paley graph; finite field; maximal clique; eigenvalue; eigenfunction; affine plane; oval
\vspace{\baselineskip}
\MSC[2010] 05E18\sep 05E30\sep 15A18
\end{keyword}
\end{frontmatter}

\section{Introduction}\label{sec0}
Let $q$ be an odd prime power, $q \equiv 1(4)$.
The \emph{Paley graph} $P(q)$ is
the Cayley graph on the additive group $\mathbb{F}_{q}^+$ of the finite
field $\mathbb{F}_{q}$ with the generating set
of all squares in
the multiplicative group $\mathbb{F}_{q}^*$.
The Paley graphs $P(q)$ are known to be strongly regular with parameters
$(q, \frac{q-1}{2}, \frac{q-5}{4}, \frac{q-1}{4})$.
The well-known Delsarte bound \cite{D73} applied to $P(q)$ says that the
cardinality of
a largest independent set (coclique) is at most $\sqrt{q}$.
Since the Paley graphs are self-complementary, the same bound holds for
a largest clique of $P(q)$.

The problem of finding clique (independence) number of Paley graphs is
open in general.
In \cite{GS16}, the Delsarte bound was improved
for infinitely many parameter tuples that correspond to Paley graphs.

Further we consider only the particular case when finite field
$\mathbb{F}_{q^2}$ is a quadratic extension over $\mathbb{F}_q$,
where $q$ is any odd prime power.
In this case the subfield $\mathbb{F}_q$ of $\mathbb{F}_{q^2}$ gives a
clique of order $q$, which meets the Delsarte bound.
In $1984$, Blokhuis \cite{B84} determined all cliques and all cocliques
of size $q$ in $P(q^2)$ and showed that they are affine images of the
subfield $\mathbb{F}_q$.



In 1996, maximal cliques of order $\frac{q+1}{2}$ and $\frac{q+3}{2}$ for $q \equiv 1(4)$ and $q \equiv 3(4)$, respectively,
were found \cite{BEHW96} by Baker et al., but an exhaustive computer search done by them
showed that these cliques are not the only cliques of such size. Moreover, there are no known maximal cliques whose size
belongs to the gap from $\frac{q+1}{2}$ (from $\frac{q+3}{2}$, respectively) to $q$.

Let $\theta$ be an eigenvalue of a graph $\Gamma$. A real-valued function on the vertex set of $\Gamma$ is called an \emph{eigenfunction}
of the graph $\Gamma$ corresponding to the eigenvalue $\theta$, if it has
at least one non-zero value and
for any vertex $\gamma$ in $\Gamma$ the condition
\begin{equation}\label{LocalCondition}
\theta\cdot f(\gamma)=\sum_{\substack{\delta\in{\Gamma(\gamma)}}}f(\delta)
\end{equation}
holds, where $\Gamma(\gamma)$ is the set of neighbours of the vertex $\gamma$.
Note that, given eigenvalue $\theta$ of a graph $\Gamma$,
a vector consisting of values of an eigenfunction of $\Gamma$ corresponding to the eigenvalue $\theta$
is an eigenvector of the
adjacency matrix of this graph corresponding to the eigenvalue $\theta$, where the values of the eigenfunction
and the indexes of the adjacency matrix have matched ordering.

There are several papers devoted to the extremal problem of studying graph eigenfunctions
with minimum cardinality of support (for more details and motivation, see \cite{KMP16}).
In \cite{V17}, Valyuzhenich found the minimum cardinality of support of an eigenfunction
corresponding to the largest non-principal eigenvalue of a Hamming graph $H(n,q)$
and characterised such eigenfunctions with the minimum cardinality of support.
In \cite{VMV}, Vorob'ev, Mogilnykh and Valyuzhenich, for all eigenvalues of a Johnson graph $J(n,\omega)$, characterised eigenfunctions
with minimum cardinality of support, where $n$ is sufficiently large.
In \cite{KMP16}, the weight-distribution lower bound for cardinality of support of an eigenfunction
of a distance-regular graph is discussed. It follows from  \cite[Corollary 1]{KMP16} that an eigenfunction of $P(q^2)$
corresponding to the eigenvalue $\theta_2 = \frac{-1-q}{2}$ has at least $q+1$ non-zero values. Since $P(q^2)$ is self-complementary,
the same bound holds for an eigenfunction of $P(q^2)$ corresponding to the eigenvalue $\theta_1 = \frac{-1+q}{2}$.


In this paper we present a new family
of maximal cliques in Paley graphs of square order.
Moreover, we use these new maximal cliques to construct a family of eigenfunctions of Paley graphs of square order,
whose cardinality of support meets the weight-distribution bound. This implies that the weight-distribution bound is tight
in the case of Paley graphs of square order.

The paper is organized as follows. In Section \ref{Preliminaries} we recall some basic notation and some preliminary results.
In Section \ref{NewMaximalCliques} we construct new maximal cliques in Paley graphs.
In Section \ref{FamilyOfEigenfunctions}, for both non-principal eigenvalues of $P(q^2)$, we present a family of eigenfunctions
with the minimum cardinality of support.

\section{Preliminaries}\label{Preliminaries}

In this section we list some useful notation and results.

\subsection{Affine plane $A(2,q)$}
Let $q$ be an odd prime power. Denote by $A(2,q)$ the point-line incidence structure, whose points are the vectors of
$2$-dimensional vector space $V(2,q)$ over $\mathbb{F}_q$,
and the lines are the additive shifts of $1$-dimensional subspaces of $V(2,q)$. It is well-known that $A(2,q)$ satisfies the axioms of
a finite affine plane of order $q$. In particular, each line contains $q$ points and there exist $q+1$ lines through a point.
An \emph{oval} in the affine plane $A(2,q)$ is a set of $q+1$ points such that no three are on a line.
A line meeting an oval in one point (in two points) is called \emph{tangent} (\emph{secant}).
For any point of an oval there exists a unique tangent at this point and $q$ secants.
By Qvist's theorem (see, for example, \cite[p. 147]{D68}), given a point that does not belong to an oval in a projective plane (and, consequently, in an affine plane) of odd order, there are either $0$ or $2$ tangents to the oval through this point.

\subsection{Finite fields of square order}
Let $d$ be a non-square in $\mathbb{F}_{q}^*$. The elements of the finite field of order $q^2$ can be considered as $\mathbb{F}_{q^2} =
\{x+y\alpha~|~x,y \in \mathbb{F}_q\}$, where $\alpha$ is a root of the polynomial $f(t) = t^2 - d$.
Since $\mathbb{F}_{q^2}$ is $2$-dimensional vector space over $\mathbb{F}_q$, we can assume that the points of $A(2,q)$ are the elements of
$\mathbb{F}_{q^2}$ and a line $l$ is presented by the elements $\{x_1+y_1\alpha + c(x_2+y_2\alpha)\}$, where $x_1+y_1\alpha \in \mathbb{F}_{q^2}$,
$x_2+y_2\alpha \in \mathbb{F}_{q^2}^*$
are fixed and $c$ runs over $\mathbb{F}_q$. The element $x_2+y_2\alpha$ is called the \emph{slope} of line $l$.
A line $l$ is called \emph{quadratic} (\emph{non-quadratic}), if its slope is a square (non-square) in $\mathbb{F}_{q^2}^*$.
Let $\beta$ be a
primitive element of the finite field $\mathbb{F}_{q}.$
Since the elements of $\mathbb{F}_{q}^* = \langle\beta^{q+1}\rangle$ are squares in $\mathbb{F}_{q^2}^*$,
the difference between any two points of quadratic (non-quadratic) line is a square (non-square) in $\mathbb{F}_{q^2}^*$.

\begin{lemma}\label{LinesThroughAPoint}
For any point of $A(2,q)$, there exists $(q+1)/2$ quadratic and $(q+1)/2$ non-quadratic lines through this point.
\end{lemma}

For any $\gamma = x+y\alpha\in \mathbb{F}_{q^2}^*$ define the \emph{norm} mapping $N$ by $N(\gamma) = \gamma^{q+1} =
\gamma\gamma^{q} = (x+y\alpha)(x-y\alpha) = x^2 - y^2d$.
The norm mapping is a homomorphism from $\mathbb{F}^*_{q^2}$ to
$\mathbb{F}^*_{q}$ with $Im(N) = \mathbb{F}^*_q$. Thus, the kernel $Ker(N)$ is the subgroup of order $q+1$
in $\mathbb{F}_{q^2}^*$. Since the kernel $Ker(N)$ is defined by the quadratic equation $x^2 - y^2d = 1$,
each line of $A(2,q)$ has at most $2$ points (elements) of $Ker(N)$. Thus, the points of $Ker(N)$ form an oval by definition.

Now we make some remarks on squares in finite fields.
\begin{lemma}\label{sq}
{\rm (1)} The element $-1$ is a square in $\mathbb{F}_q^*$ iff $q \equiv 1(4)$;\\
{\rm (2)} For any non-square $d$ in $\mathbb{F}_q^*$ the element $-d$ is a square in $\mathbb{F}_q^*$ iff $q \equiv 3(4)$.
\end{lemma}

The following lemma can be used to test whether an element $\gamma = x+y\alpha \in \mathbb{F}_{q^2}^*$ is a square.

\begin{lemma}[\cite{BEHW96}, Lemma 2]\label{sqq}
 An element $\gamma = x+y\alpha \in \mathbb{F}_{q^2}^*$ is a square iff $N(\gamma) = x^2-y^2d$ is a square in $\mathbb{F}_q^*$.
\end{lemma}

Lemma \ref{sqqq} immediately follows from Lemma \ref{sqq}, Lemma \ref{sq} and the fact that $N(\alpha) = -d$.
\begin{lemma}\label{sqqq}
The element $\alpha$ is a square in $\mathbb{F}_{q^2}^*$ iff $q \equiv 3(4)$.
\end{lemma}

\section{Family of new maximal cliques}\label{NewMaximalCliques}
In this section, for any odd prime power $q$, we give a construction of maximal cliques in the Paley graph $P(q^2)$.

Let $\beta$ be a primitive element of the finite field $\mathbb{F}_{q^2}$.
Put $\omega:=\beta^{q-1}$. Note that $\omega$ is a square in $\mathbb{F}_{q^2}^*$. Then
the subgroup of order $q+1$ in $\mathbb{F}_{q^2}^*$ is presented by $Q = \langle\omega\rangle$.
Put $Q_0 := \{1, \omega^2, \omega^4, \ldots, \omega^{q-1}\}, Q_1 := \{\omega, \omega^3, \omega^5 \ldots, \omega^q\}.$
We have $Q = Q_0 \cup Q_1$. It was noticed in Section \ref{Preliminaries} that the elements of $Q$ form an oval as
points of $A(2,q)$.

The main goal of this section is to prove the following result.

\begin{theorem}\label{main1}~\\
{\rm(1)} If $q \equiv 1 (4)$, then $Q_0$ and $Q_1$ are maximal cocliques of size $\frac{q+1}{2}$ in the graph $P(q^2)$;\\
{\rm(2)} If $q \equiv 3 (4)$, then $Q_0 \cup \{0\}$ and $Q_1 \cup \{0\}$ are maximal cliques of size $\frac{q+3}{2}$ in the graph $P(q^2)$.
\end{theorem}

It follows from Theorem \ref{main1} that, for any $q \equiv 1 (4)$, the sets $Q_0$ and $Q_1$ induce maximal cliques of size $\frac{q+1}{2}$ in the complementary graph $\overline{P(q^2)}$, which is isomorphic to $P(q^2)$.

For any $\beta_1 \in \mathbb{F}_{q^2}^*$, $\beta_2 \in \mathbb{F}_{q^2}$ define an affine transformation
$\psi_{\beta_1, \beta_2}  :\mathbb{F}_{q^2} \rightarrow \mathbb{F}_{q^2}$ by the following rule
$$\psi_{\beta_1, \beta_2}(\gamma) = \beta_1\gamma + \beta_2.$$
The set $T := \{\psi_{\beta_1,\beta_2}~|~\beta_1 \in (\mathbb{F}_{q^2}^*)^2\}$ forms a subgroup in ${\rm Aut}(P(q^2))$.

Let us consider the sets $T_Q:=\{\psi_{\beta_1,0}~|~\beta_1 \in Q\} = \langle\psi_{\omega,0}\rangle,$
$T_{Q_0}:=\{\psi_{\beta_{1},0}~|~\beta_1 \in Q_0\} = \langle\psi_{\omega^2,0}\rangle, $
$T_{Q_1}:=\{\psi_{\beta_{1},0}~|~\beta_1 \in Q_1\},$ where $T_{Q_0}$ and $T_Q$ are subgroups of $T$
such that $T_{Q_0} < T_Q < T$ holds.

\begin{lemma}\label{T_Q}
The following statements hold.\\
{\rm(1)} $T_Q$ preserves the lines of $A(2,q)$;\\
{\rm(2)} $T_Q$ stabilizes setwise the oval $Q$;\\
{\rm(3)} $T_Q$ acts transitively on the points of $Q$;\\
{\rm(4)} $T_Q$ acts transitively on the tangents to the oval $Q$;\\
{\rm(5)} $T_{Q_0}$ stabilizes setwise the sets $Q_0$ and $Q_1$ and acts transitively on each of them;\\
{\rm(6)} each element of $T_{Q_1}$ swaps the sets $Q_0$ and $Q_1$;\\
{\rm(7)} $T_{Q}$ stabilizes setwise the sets of quadratic and non-quadratic lines of $A(2,q)$.
\end{lemma}

\begin{lemma}\label{TangentsAreTheSame}
The following statements hold.\\
{\rm(1)} If $q \equiv 1(4),$ then all tangents to the oval $Q$ are non-quadratic;\\
{\rm(2)} If $q \equiv 3(4),$ then all tangents to the oval $Q$ are quadratic.
\end{lemma}
\proof
Note that the line $l = \{1+c\alpha~|~c\in \mathbb{F}_q\}$ is a tangent to the oval $Q$ at the point $(1,0)$.
By Lemma \ref{T_Q}(4) and \ref{T_Q}(7), $T_Q$ acts transitively on the tangents to $Q$ and stabilizes the set of quadratic and non-quadratic lines.
This means that a tangent to $Q$ is quadratic iff the tangent $l$ is quadratic. The tangent $l$ has the slope $\alpha$, which is a square
iff $q \equiv 3(4)$ by Lemma \ref{sqqq}. The lemma is proved.

\begin{lemma}\label{NeighboursOf1}
The following statements hold.\\
{\rm(1)} If $q \equiv 1(4),$ then the identity $1$ is adjacent to each element of $Q_1$ and has no neighbours in $Q_0\setminus\{1\}$;\\
{\rm(2)} If $q \equiv 3(4),$ then the identity $1$ is adjacent to each element of $Q_0\setminus\{1\}$ and has no neighbours in $Q_1$.
\end{lemma}
\proof
Pick an arbitrary element $\gamma = x+y\alpha \in Q\setminus\{1\}$. Then the element $\gamma^2 = x^2+dy^2 + 2xy\alpha$
represents an arbitrary element of $Q_0\setminus\{1\}$. Moreover, since $\gamma$ is a point of the oval $Q$,
the equality $x^2-dy^2 = 1$ holds.

Let us consider the difference
$$\gamma^2-1 = x^2-1+dy^2 + 2xy\alpha = 2dy^2+2xy\alpha.$$
By Lemma \ref{sqq}, this difference is a square in $\mathbb{F}_{q^2}^*$ iff the element $N(\gamma^2-1) = N(2dy^2+2xy\alpha) =
4d^2y^4-4x^2y^2d = 4dy^2(dy^2 - x^2) = -4dy^2
$
is a square in $\mathbb{F}_q^*$. By Lemma \ref{sq}(1), the element $-4dy^2$ is a square in $\mathbb{F}_q^*$ iff $q \equiv 3(4)$.

Thus, one of the following two cases holds. If $q \equiv 1(4)$, then a non-quadratic line through $1$ either contains a point from
$Q_0\setminus\{1\}$ or
is a tangent. Consequently, a quadratic line through $1$ contains a point from $Q_1$.
If $q \equiv 3(4)$, then a quadratic line through $1$ either contains a point from
$Q_0\setminus\{1\}$ or
is a tangent. Consequently, a non-quadratic line through $1$ contains a point from $Q_1$.
The lemma is proved.

\begin{lemma}\label{AdjStructureOfQ}
The following statements hold.\\
{\rm(1)} If $q \equiv 1(4),$ then the graph induced by $Q$ is a complete bipartite with parts $Q_0$ and $Q_1$;\\
{\rm(2)} If $q \equiv 3(4),$ then the graph induced by $Q$ is a disjoint union of cliques $Q_0$ and $Q_1$.
\end{lemma}
\proof
The lemma immediately follows from Lemmas \ref{NeighboursOf1} and \ref{T_Q}(5).

\medskip
Now let us complete the proof of Theorem \ref{main1}.

(1) By Lemma \ref{T_Q}(6), $Q_0$ and $Q_1$ induce isomorphic subgraphs, so,
it is enough to prove that coclique $Q_0$ is maximal.

Suppose that the coclique $Q_0$ can be extended by an element
$\delta \in \mathbb{F}_{q^2}^* \setminus Q_0$.
Each quadratic line through the point $0$ is a secant of the oval $Q$;
such a line contains a point from $Q_0$ and a point from $Q_1$.
Thus, the line, which contains the points $0$ and $\delta$, is non-quadratic.

Consider a line $L$ that connects the point $\delta$ with a point from $Q$.
Suppose $L$ is a quadratic. Then $L$ is a secant that contains a point
from $Q_0$, which is a contradiction. So, $L$ is non-quadratic.

Consider all lines connecting $\delta$ with each point of the oval $Q$.
Since each such a line is non-quadratic, we obtain $(q +1)/2$ secants
of $Q$; these $(q + 1)/2$ secants are exactly $(q + 1)/2$ non-quadratic
lines through the point $\delta$; one of them contains the point $0$.
Since each secant of $Q$ through the point $0$ is quadratic, we obtain a
contradiction.

(2)
Recall that $Q$ is the set of all elements with the norm $1$.
For any element $s \in \mathbb{F}_q^*$
the set $sQ = sQ_0 \cup sQ_1$ is the set of all elements with the
quadratic norm $s^2$;
the sets $sQ_0 \cup \{s\}$ and $sQ_1 \cup \{s\}$ induce cliques of order
($q + 3)/2$.
Each square of $\mathbb{F}_{q^2}^*$ belongs to a unique set $sQ$ for some $s \in \mathbb{F}_q^*$.
Suppose that the clique $Q_0 \cup \{0\}$ can be extended by a square
$\delta \in sQ_i$ , where $i \in \{0, 1\}$ and $s \in \mathbb{F}_q^*, s \ne 1$.
Since $T_{Q_0}$ stabilizes (setwise) the clique $Q_0 \cup \{0\}$ and acts
transitively on the set of vertices of the clique $sQ_i$, each vertex of
the clique $sQ_i$ is adjacent to each vertex of the clique $Q_0 \cup \{0\}$.
We obtain the clique $sQ_i \cup Q_0 \cup \{0\}$ of the size
$|sQ_i \cup Q_0 \cup \{0\}| = |sQ_i|+|Q_0|+1 = (q+1)/2+(q+1)/2+1 = q+2 > q$,
which is a contradiction. The theorem is proved.

\medskip
An exhaustive computer search shows that the maximal cliques presented in \cite{BEHW96} and
the maximal cliques given by Theorem \ref{main1} are not the only maximal cliques of size
$\frac{q+1}{2}$ ($\frac{q+3}{2}$, respectively) in
$P(q^2)$.
\section{Family of eigenfunctions with minimum cardinality of support}\label{FamilyOfEigenfunctions}

The goal of this section is to prove the following theorem.
\begin{theorem}
 Let  $f:\mathbb{F}_{q^2} \longrightarrow{\mathbb{R}}$ be a function defined by the following rule.
$$
f(\gamma):=\begin{cases}
1,&\text{if $\gamma \in Q_0$;}\\
-1,&\text{if $\gamma \in Q_1$;}\\
0,&\text{otherwise.}
\end{cases}
.$$
{\rm(1)} If $q \equiv 1 (4)$, then $f$ is an eigenfunction of $P(q^2)$ corresponding to the eigenvalue $\theta_2 = \frac{-1-q}{2}$ and $|Supp(f)| = q+1$ holds;\\
{\rm(2)} If $q \equiv 3 (4)$, then $f$ is an eigenfunction of $P(q^2)$ corresponding to the eigenvalue $\theta_1 = \frac{-1+q}{2}$ and $|Supp(f)| = q+1$ holds.
\end{theorem}
\proof
We prove the theorem by checking the equality (\ref{LocalCondition}) for each vertex of $P(q^2)$.\\
(1) Let $\gamma$ be a vertex of $P(q^2)$ that represents a point of the oval $Q$.
By Lemma \ref{AdjStructureOfQ}(1), the vertices of $Q$ induce a complete bipartite graph with parts $Q_0$ and $Q_1$
of cardinality $\frac{1+q}{2}$. Assume that $\gamma$ belongs to $Q_i$, where $i\in\{0,1\}$.
Then the neighbours of $\gamma$ with non-zero value of $f$ are exactly elements of $Q_{1-i}$.
So, for any vertex $\gamma \in Q$ the equality $(\ref{LocalCondition})$ holds.

Suppose $\gamma$ does not belong to the oval $Q$. Note that $f(\gamma) = 0$ holds.
We consider the set $L_\gamma$ of lines that connect the point $\gamma$ with a point of $Q$.
Then $L_\gamma$ contains either $0$ or $2$ tangents to the oval $Q$.

If $L_\gamma$ contains no tangents to the oval $Q$,
then each line in $L_\gamma$ is a secant of the oval $Q$. So, the neighbours of $\gamma$ with non-zero value of $f$ are exactly
those points of $Q$ that belong to quadratic lines of $L_\gamma$. Each such quadratic line contains
some pair of adjacent vertices of $Q$.
By Lemma \ref{AdjStructureOfQ}, any edge in $Q$ connects a vertex from $Q_0$ and a vertex from $Q_1$.
Thus, the equality $(\ref{LocalCondition})$ holds.

If $L_\gamma$ contains two tangents to the oval $Q$, then by Lemma \ref{TangentsAreTheSame}(1) these tangents are non-quadratic
and, indeed, $\gamma$ is not adjacent in the graph $P(q^2)$ to the two points of tangency.
Each line in $L_\gamma$ excepting the two tangents is either a secant of the oval $Q$ or has no points of $Q$.
So, the neighbours of $\gamma$ with non-zero value of $f$ are exactly
those points of $Q$ that belong to quadratic secant $L_\gamma$. Each such quadratic secant contains
some pair of adjacent vertices of $Q$.
By Lemma \ref{AdjStructureOfQ}, any edge in $Q$ connects a vertex from $Q_0$ and a vertex from $Q_1$.
Thus, the equality $(\ref{LocalCondition})$ holds.

(2) Let $\gamma$ be a vertex of $P(q^2)$ that represents a point of the oval $Q$.
By Lemma \ref{AdjStructureOfQ}(2), the vertices of $Q$ induce the disjoint pair of cliques $Q_0$ and $Q_1$
of cardinality $\frac{1+q}{2}$. Assume that $\gamma$ belongs to $Q_i$, where $i\in\{0,1\}$.
Then the neighbours of $\gamma$ with non-zero value of $f$ are exactly elements of $Q_{i} \setminus \{\gamma\}$.
So, for any vertex $\gamma \in Q$ the equality $(\ref{LocalCondition})$ holds.

Suppose $\gamma$ does not belong to the oval $Q$. Note that $f(\gamma) = 0$ holds.
We consider the set $L_\gamma$ of lines that connect the point $\gamma$ with a point of $Q$.
Then $L_\gamma$ contains either $0$ or $2$ tangents to the oval $Q$.

If $L_\gamma$ contains no tangents to the oval $Q$,
then each line in $L_\gamma$ is either a secant of the oval $Q$ or has no points of $Q$.
So, the neighbours of $\gamma$ with non-zero value of $f$ are exactly
those points of $Q$ that belong to quadratic secants of $L_\gamma$. Each such quadratic secant contains
a unique pair of adjacent vertices of $Q$.
By Lemma \ref{AdjStructureOfQ}, any edge in $Q$ connects a pair of vertices either both from $Q_0$ or both from $Q_1$.
This means that $\gamma$ is adjacent to each vertex of $Q$.
Thus, the equality $(\ref{LocalCondition})$ holds.

If $L\gamma$ contains two tangents to the oval $Q$,
then by Lemma \ref{TangentsAreTheSame}(2) these tangents are quadratic.
Each line in $L\gamma$ excepting the two tangents is either a secant of the oval $Q$ or has no points of $Q$.
So, the neighbours of $\gamma$ with non-zero value of $f$ are exactly
those points of $Q$ that belong to quadratic secants of $L_\gamma$.
Any non-quadratic line from $L_\gamma$ is a secant of the oval $Q$ and connects a vertex from $Q_0$ with a vertex $Q_1$.
This means that non-quadratic lines from $L_\gamma$ contain the same number of points from $Q_0$ and $Q_1$.
Since $|Q_0| = |Q_1|$, quadratic lines from $L_\gamma$ contain the same number of points from $Q_0$ and $Q_1$.
Thus, the equality $(\ref{LocalCondition})$ holds. The theorem is proved.

\section*{Acknowledgment} \label{Ack}
We thank Alexander Gavrilyuk, Elena Konstantinova and Denis Krotov for useful suggestions,
which significantly improved this paper.

\end{document}